\documentclass{article}
\usepackage{color}
\usepackage{url}
\usepackage{graphicx}
\usepackage{amsmath,amssymb}
\usepackage{epsfig,graphicx,epstopdf}
\allowdisplaybreaks
\usepackage{float}
\usepackage{graphicx}
\usepackage{appendix}
\usepackage{listings}
\usepackage{subfigure}
\usepackage{cite}
\usepackage{bm}
\usepackage[numbers]{natbib}
\usepackage{hyperref}
\usepackage{mathrsfs}

\usepackage{makecell}
\usepackage{multirow}
\usepackage{algorithm}
\usepackage{algorithmicx}
\usepackage{algpseudocode}
\usepackage[T1]{fontenc}
\usepackage{authblk}

\usepackage{multicol}

\newtheorem{theorem}{Theorem}[section]
\newtheorem{lemma}{Lemma}[section]
\newtheorem{remark}{Remark}[section]

\newtheorem{example}{Example}[section]
\numberwithin{equation}{section}

\newenvironment{proof}{\paragraph{Proof:}}{\hfill$\square$}
\makeatletter

\title{A modified dynamic diffusion finite element method with optimal convergence rate for convection-diffusion-reaction equations\thanks{This work was supported in part by the Scientific and Technological Research Program of Chongqing Municipal Education Commission (KJZD-M202300705), the National Natural Science Foundation
of China (12171340), and the Graduate Research and Innovation Fund Project of Chongqing Jiaotong University (2023S0107).}}

\author{Shaohong Du\thanks{
School of Mathematics and Statistics, Chongqing Jiaotong University, Chongqing 400074, China. E-mail: {duzheyan.student@sina.com}.}, \quad
Qianqian Hou\thanks{School of Mathematics and Statistics, Chongqing Jiaotong University, Chongqing 400074, China. E-mail: {1719792019.com}.}, \quad
Xiaoping Xie\thanks{Corresponding author. School of Mathematics, Sichuan University, Chengdu 610064, China.  Email: {xpxie@scu.edu.cn}}
}

\date{}
\begin{document}

\maketitle
\renewcommand{\thefootnote}{\fnsymbol{footnote}}

\begin{abstract}
In this paper, we develop a modified nonlinear dynamic diffusion (DD) finite element method for convection-diffusion-reaction equations. This method is free of stabilization
parameters and is capable of precluding spurious oscillations. We prove   existence  and,  under an assumption of small mesh size,   uniqueness of the discrete solution, and derive  the  optimal first oder convergence rate of  the approximation error in  the energy norm   plus a dissipative term. Numerical examples are provided to verify the theoretical analysis.
\end{abstract}

{\it Keywords}: convection-diffusion-reaction equation, dynamical diffusion method, optimal convergence rate.

{\it AMS subject classifications}: 65K10, 65N30, 65N21, 49M25, 49K20.
\section {Introduction}
Let $\Omega\subset\mathbb{R}^{d} $ ($d=2,3$) be a bounded domain with Lipschitz boundary $\Gamma=\bar{\Gamma}_{D}\cup\bar{\Gamma}_{N}$, $\Gamma_{D}\cap\Gamma_{N}=\emptyset$. We consider the following   convection-diffusion-reaction model:
\begin{equation}\label{formula1}
\left \{\begin{array}{ll}
-\varepsilon\triangle u+{\beta}\cdot\nabla u+\sigma u = f\ \ &{\rm in}\ \Omega,\vspace{2mm}\\
u= 0\ \ {\rm on}\ \Gamma_{D},\ \varepsilon\nabla u\cdot{\bf n}=g\ \ &{\rm on}\ \Gamma_{N},
\end{array}\right.
\end{equation}
where  $u$ represents the quantity being transported, $0<\varepsilon\ll1$ the (constant) diffusivity, $\beta\in[W^{1,\infty}(\Omega)]^{d}$ the velocity
field, $0\leq\sigma\in L^{\infty}(\Omega)$ the reaction coefficient, $f\in L^{2}(\Omega)$ the source term, ${\bf n}$ the outward unit normal vector along $\Gamma$, $g\in L^{2}(\Gamma_{N})$ the Neumann boundary condition.

We first make the following two assumptions on the model problem:
\begin{itemize}
\item[(A1)]\ The Dirichlet boundary $\Gamma_{D}$ has a positive $(d-1)-$dimensional Lebesgue measure and includes the inflow boundary, i.e.  $\{{\bf x}\in\Gamma : \beta({\bf x})\cdot{\bf n}<0\}\subset\Gamma_{D}$;
\item[(A2)]\ There are two nonnegative constants $\gamma$ and $c_{\sigma}$, independent of $\varepsilon$, such that
\begin{equation}\label{formula2}
\sigma-\frac12 \nabla\cdot\beta\geq\gamma\ \ {\rm and}\ \ \|\sigma\|_{0,\infty}\leq c_{\sigma}\gamma.
\end{equation}
\end{itemize}
Assumption (A2)  ensures that we can simultaneously handle the cases of $\sigma=0$ (cf. \cite{Verfurth2005}) and $\sigma\neq 0$, with the formmer case corresponding to $\gamma=0$ and $c_{\sigma}=0$.

It is known that standard finite element methods (FEMs) usually give rise to  nonphysical solutions with numerical oscillations for the
convection-dominated situation of \eqref{formula1}
due to local singularities arising from interior or boundary layers  \cite{Chalot1998,Carey1983,Roos1996,Roos2003,Franca2006,Cangiani2004}.
 To alleviate the effects of  numerical oscillations, there have developed many linearly stabilized methods
 such as the Streamline-Upwind$/$Petrov-Galerkin (SUPG) method \cite{Brooks1982}, the Galerkin-Least-Squares method \cite{HUghes1989},   the unusual stabilized finite element method \cite{Franca},  
  discontinuous Galerkin (DG) methods  \cite{b2Cockburn, b1Ayuso, bs4Houston, b5Zarin,  b7Buffa}, hybridizable DG methods \cite{b8Chen,b9Chen,bc,b10Fu}, and weak Galerkin methods \cite{MWG, CFX,WG1}.  We mention that for these linearly stabilized
methods,  the spurious oscillations may still remain in the neighbourhood of sharp gradients and a suitable design of  the stabilization parameters is required.

The Variational Multiscale (VMS) approach \cite{Hughes1995,Hughes2004,Hughes;1998,Hughes2000,Layton2002} provides a procedure for deriving models and
numerical methods capable of dealing with multiscale phenomena.
The basic idea of VMS is to decompose the variable of interest into a resolved coarse scale and an unresolved subgrid scale and, in so doing, the subgrid variability can be captured in the numerical solution. Many stabilized methods can be incorporated into  the VMS framework, e.g. the residual-free bubble method \cite{Brezzi1994,Franca2006,Barreda2019}, the multiscale finite element method \cite{Hou1997,Coutinho;2011,Girault1986,Gravemier2003}, the
subgrid stabilization method \cite{Brezzi2000,Guermond1999,Guermond2001,Hauke2001,John2006,Guermond2004} and the nonlinear
subgrid method \cite{Knopp2002,Arruda2010,Du2023,Santos2007,Santos2012,Santos2018,Santos2021}. In particular, by introducing a nonlinear artificial dissipation dynamically adjusted in terms of the notion of scale separation the nonlinear subgrid method, also called the Dynamic Diffusion (DD) method, is able to preclude local and
nonlocal oscillations in the convection-diffusion problems without  requiring any
tuned-up stabilization parameter.  It should be pointed out that  the existence of approximate solution is not guaranteed  for the  DD scheme in \cite{Santos2007},
as the local nonlinear operator defined by the artificial diffusion does not have an upper bound, and that an improvement of the artificial diffusion was given in \cite{Santos2021} which allows for the upper and lower bounds of the corresponding local nonlinear operator.

Since the nonlinear two-scale DD method does not have some important properties like the Lipschitz continuity, the linearity preservation and the discrete maximum principle  \cite{Burman2005}, its numerical analysis shares
the difficulty of classical nonlinear shock-capturing methods \cite{condina1993,Galeao1988,Hughes1986,John2007}: the first is the uniqueness of approximation solution; and the second is that only a convergence
rate of order $O(h^{1/2})$  can be proven   in the energy norm for the linear element approximation,   though   numerical experiments demonstrating  the optimal convergence rate $O(h)$.

In this paper, we develop a modified nonlinear DD finite element method for the convection-diffusion-reaction model \eqref{formula1}, give results of  existence and uniqueness for discrete solution, and derive a priori error estimates. Compared with \cite{Santos2021},
our main contributions lie in the following  aspects:

\begin{itemize}
\item The considered model \eqref{formula1} is more general where the velocity $\beta$ is not necessarily divergence free and the reaction coefficient $\sigma$ is a function, while
in \cite{Santos2021} it has been assumed that $\nabla\cdot\beta=0$ and $\sigma$ is a constant;

\item We modify the artificial diffusivity term to a proper dimensional scale, thus unify the two alternative variants of the artificial
diffusion in \cite{Santos2021} into one form. In the proof of the existence of approximation solutions, we do not assume the boundedness of the maximum norm of the gradient of functions in approximation space (see Remark 4.1). We further prove the uniqueness of approximation solution for the modified DD method under an assumption on small mesh size.

\item We prove the optimal convergence rate $O(h)$ in the energy norm of approximation error plus a dissipative term for the modified DD method.
\end{itemize}

This paper is organized as follows. Section 2 gives notations and weak formulations of the model problem. Section 3 presents the modified DD scheme. Section 4 proves the existence and uniqueness of discrete solutions for the modified DD method, and also contains some
preliminary results at the beginning of this section. Section 5 analyzes the optimal
convergence order for the modified DD method. Finally, section 6 provides some numerical tests to support our theoretical results.

\section{Notation and weak formulation}
For any subdomain $\omega$ of $\Omega$, denote $\|\cdot\|_{m,p,\omega}$ (resp. $|\cdot|_{m,p,\omega}$) the standard norm (resp. semi-norm) in the Sobolev space $W^{m,p}(\omega)$. In
case of $p=2$, $W^{m,2}(\omega)=H^{m}(\omega)$ ($H^{0}(\omega)=L^{2}(\omega)$) is the Hilbert space (with the inner product $(\cdot,\cdot)_{\omega}$) and is
equipped with norm $\|\cdot\|_{m,\omega}$ (resp. $|\cdot|_{m,\omega}$ ); in case of $m=0,p=\infty$, $W^{0,\infty}(\omega)=L^{\infty}(\omega)$ is the Lebesgue space with norm $\|\cdot\|_{0,\infty,\omega}$. We
omit the symbol $\Omega$ in the notations above if $\omega = \Omega$, and also set
\begin{equation*}
H_{D}^{1}(\Omega)=\{v\in H^{1}(\Omega): v|_{\Gamma_{D}}=0\}.
\end{equation*}

Let $\mathcal{T}_{h}$ be a shape regular triangulation of $\Omega$ into triangles/tetrahedrons, and denote by $h$ the mesh-size function with $h|_{T}=h_{T} :={\rm diam}(T)$ for
any $T\in \mathcal{T}_{h}$. Let $V_{h}\subset H^{1}(\Omega)$ be the usual linear element space, and $\psi_{T}$ be the bubble function of degree $d+1$ defined in each element $T$. Introduce finite dimensional spaces
   $$V_{b}:=\oplus|_{T\in \mathcal{T}_{h}}{\rm span}\{\psi_{T}\}, \quad  V_{h}^{D}:=V_{h}\cap H_{D}^{1}(\Omega), \quad   V_{hb}^{D}:=V_{h}^{D}\oplus V_{b}.$$
Notice that in a two-scale method $V_{b}$ is the unresolved (fine) scale space, $V_{h}^{D}$ is the resolved (coarse) scale space, whereas $V_{b}$ is defined by using the bubble functions due to reduction of the computational cost.

Throughout this paper, $c_{i},i=0,1,2,\cdots$, stand for generic positive constants independent of the (constant) diffusivity $\varepsilon$ and mesh size $h$, and may take different values at different occurrences.

We introduce a bilinear form $B(\cdot,\cdot)$ on $H_{D}^{1}(\Omega)\times H_{D}^{1}(\Omega)$ by
\begin{equation}\label{Buv}
B(u,v):=\varepsilon(\nabla u,\nabla v)+({\beta}\cdot\nabla u, v)+(\sigma u,v),\ \forall u,v\in H_{D}^{1}(\Omega).
\end{equation}
Then the standard variational formulation of \eqref{formula1} reads: find $u\in H_{D}^{1}(\Omega)$ such that
\begin{equation}\label{formula3}
B(u,v)=(f,v)+(g,v)_{\Gamma_{N}}, \ \ \forall v\in H_{D}^{1}(\Omega).
\end{equation}

For any $v\in H_{D}^{1}(\Omega)$, set
\begin{equation}
\label{|||-T}
|||v|||_{T}  :=\left(\varepsilon|v|_{1,T}^{2}+\gamma\|v\|_{0,T}^{2}\right)^{1/2}, \quad \forall T\in \mathcal{T}_{h}\end{equation}
and define an energy norm $|||\cdot|||$ by
\begin{equation}
\label{|||-norm}
|||v|||^{2} :=\sum\limits_{T\in \mathcal{T}_{h}} |||v|||_{T}^2=\varepsilon|v|_{1}^{2}+\gamma\|v\|_{0}^{2}.
\end{equation}
From (A1) and (A2) we easily obtain the coercivity and boundedness of $B(\cdot,\cdot)$, i.e. for $\forall v,w\in H_{D}^{1}(\Omega)$,
\begin{equation}\label{formula4}
B(v,v)\geq|||v|||^{2}
\end{equation}
and
\begin{equation}\label{formula5}
|B(v,w)|\leq M_1|||v|||\ |||w|||
\end{equation}
with $$M_1:=\left\{\begin{array}{ll}
1+\|\sigma\|_{0,\infty}+\frac{\|\beta\|_{0,\infty}}{\sqrt{\varepsilon\gamma}} , & \text{if } \gamma>0,\\
1+ \frac{c_0\|\beta\|_{0,\infty}}{ \varepsilon } , &\text{if } \gamma=0,
\end{array}
\right.$$
where we have used the fact
\begin{equation}\label{formula0-1}
||v||_0\leq c_0|v|_1, \quad  v \in H_{D}^{1}(\Omega).
\end{equation}
 These two conditions ensure that (\ref{formula3}) admits a unique solution $u\in H_{D}^{1}(\Omega)$.

\section{A modified DD method}

Define a linear projection operator $\kappa_{h} :V_{hb}^{D}\rightarrow V_{h}^{D}$ by
\begin{equation}\label{december1}
\kappa_{h}(w_{hb})=w_{h}, \quad \forall w_{hb}=w_{h}+w_{b} \text{ with }w_{h}\in V_{h}^{D},w_{b}\in V_{b}
\end{equation}
and introduce, for any $T\in\mathcal{T}_{h}$, a local trilinear operator  $D_{T}(\cdot;\cdot,\cdot): V_{hb}^{D}\times V_{hb}^{D} \times V_{hb}^{D}\rightarrow\mathbb{R}$ given by
\begin{equation}\label{formulatwo-}
D_{T}(w_{hb};u_{hb},v_{hb}):=\displaystyle\int_{T}\xi_{T}(\kappa_{h}(w_{hb}))\nabla u_{hb}\cdot\nabla v_{hb}d{\bf x},
\end{equation}
where $\xi_T: V_{h}^D\rightarrow \mathbb{R}^{+}$ denotes an additional diffusivity  coefficient to be determined later.

The DD method for the convection-diffusion-reaction model \eqref{formula1} reads \cite{Santos2018}: Find $u_{hb}=u_{h}+u_{b}\in V_{hb}^{D}$ with $ u_{h}\in V_{h}^{D}, u_{b}\in V_{b}$, such that
\begin{equation}\label{formulatwo0}
\displaystyle B(u_{hb},v_{hb})+\sum\limits_{T\in\mathcal{T}_{h}}D_{T}(u_{hb};u_{hb},v_{hb})=(f,v_{hb})+(g,v_{hb})_{\Gamma_{N}},\ \forall v_{hb}\in V_{hb}^{D},
\end{equation}
Note that the term $D_{T}(\cdot;\cdot,\cdot)$ introduces an isotropic diffusion on all scales, whose magnitude is dynamically determined by
imposing the restriction $\xi_T$ on the resolved scale solution $u_h$, i.e.
$$\xi_T(\kappa_{h}(u_{hb}))=\xi_T(u_{h}).$$
How to design $\xi_T(u_{h})$  in each element $T$ is  the key to the DD method.

In  \cite{Santos2018} the diffusivity restriction is defined as
\begin{equation}\label{december2}
\xi_{T}(u_{h})=\frac{1}{2}\hbar\|\beta_b\|_{0,T},
\end{equation}
where $\hbar$ is a characteristic length scale, e.g. $\hbar=\sqrt{|T|}$,
 and the subgrid velocity field $\beta_{b}$ solves the optimization problem
\begin{equation}\label{formulaone}
  \left \{\begin{array}{ll}
   {\min}\ E_{k}=\frac{1}{2}\|\beta_{b}\|_{0,T}^{2}\\
   {\rm subject\ to}\\
   -\varepsilon\triangle_h u_{h}+(\beta-\beta_{b})\cdot\nabla u_{h}+\sigma u_{h}=f\ \ {\rm in}\ T\in\mathcal{T}_{h},
  \end{array}\right.
  \end{equation}
It is easy to see that
\begin{equation}\label{december3}
\beta_{b}=\left \{\begin{array}{ll}
   \frac{R(u_{h})}{\|\nabla u_{h}\|_{0,T}^{2}}\nabla u_{h}\ \ {\rm if}\ \|\nabla u_{h}\|_{0,T}\neq 0,\\
   {\bf 0},\ \ {\rm otherwise},
  \end{array}\right.
  \end{equation}
where
\begin{equation}
\label{residual}
R(u_{h}):=-\varepsilon\triangle_h u_{h}+\beta\cdot\nabla u_{h}+\sigma u_{h}-f
\end{equation}
 is the residual of the resolved scale solution. In particular, it holds
\begin{align} \label{Ruh}
R(u_{h})= \beta\cdot\nabla u_{h}+\sigma u_{h}-f
\end{align}
since $\triangle_h u_{h}=0$ for $u_{h}\in V_{h}^{D}$. We note that the residual term $R(u_{h})$ was replaced  in  \cite{Santos2018} by
    $|R(u_{h})|$ which   was actually evaluated at the barycenter of each element $T\in\mathcal{T}_{h}$.

 Since the local nonlinear operator $D_{T}(\cdot;\cdot,\cdot)$ defined by the artificial diffusion $\xi_{T}(\cdot)$ is unbounded
 due to the fractor  $\frac{1}{\|\nabla u_{h}\|_{0,T}^2}$ in \eqref{december3}, the existence of approximate solutions to the DD scheme \eqref{formulatwo}
 is not guaranteed, let alone the analysis of stability and convergence. A remedy was proposed in \cite{Santos2021} to choose
\begin{equation}\label{formula6-}
\xi_{T}(u_{h})=\left \{\begin{array}{ll}
   h_{T}\frac{\|R(u_{h})\|_{0,T}}{\|u_{h}\|_{1,T}+\tau}\ \ {\rm if}\ P_{e_{T}}>1,\\
    0,\ \ {\rm otherwise}
  \end{array}\right. \quad \text{ for }\sigma>0
\end{equation}
and
\begin{equation}\label{formula7-}
\xi_{T}(u_{h})=\left \{\begin{array}{ll}
   h_{T}\frac{\|R(u_{h})\|_{0,T}}{|u_{h}|_{1,T}+\tau}\ \ {\rm if}\ P_{e_{T}}>1 ,\\
    0,\ \ {\rm otherwise}
  \end{array}\right. \quad \text{ for }\sigma=0,
\end{equation}
where $P_{e_{T}}:=\frac{\|\beta\|_{0,T}h_{T}}{2\varepsilon}$
denotes the Pecl\'{e}t number, and $\tau$ is an appropriate constant.

In what follows we shall modify  the artificial diffusion  $\xi_{T}$ in the scheme \eqref{formulatwo0}  so as to  obtain the following  modified DD method:

Find $u_{hb}=u_{h}+u_{b}\in V_{hb}^{D}$ with $ u_{h}\in V_{h}^{D}, u_{b}\in V_{b}$, such that
\begin{equation}\label{formulatwo}
\displaystyle a(u_{hb};u_{hb},v_{hb})=(f,v_{hb})+(g,v_{hb})_{\Gamma_{N}},\ \forall v_{hb}\in V_{hb}^{D},
\end{equation}
where $a(\cdot;\cdot,\cdot): V_{hb}^{D}\times V_{hb}^{D} \times V_{hb}^{D}\rightarrow \mathbb{R}$ is defined as
$$a(w_{hb};u_{hb},v_{hb}):=B( u_{hb},v_{hb})+\sum\limits_{T\in\mathcal{T}_{h}}D_{T}(w_{hb};u_{hb},v_{hb} ) , \quad  \forall  w_{hb},u_{hb},v_{hb} \in V_{hb}^{D},$$
and  the   artificial diffusion $\xi_T(\kappa_{h}(w_{hb}))=\xi_T(w_{h})$  in the definition \eqref{formulatwo-} of $D_{T}(\cdot;\cdot,\cdot)$  is given by
\begin{equation}\label{formula6}
\xi_{T}(w_{h}):=\left \{\begin{array}{ll}
   \frac{h_{T}\|R(w_{h})\|_{0,T}}{A_T(w_{h})+\tau}\ & {\rm if}\ P_{e_{T}}>1,\\
    0,\ & {\rm otherwise}
  \end{array}\right.
\end{equation}
with
\begin{equation}\label{choice-A}
 A_T(w_{h}):= \|\beta\|_{0,\infty,T}|w_{h}|_{1,T}+\|\sigma\|_{0,\infty,T}\|w_{h}\|_{0,T}
\end{equation}
and
\begin{equation}\label{choice-tau}
\tau:=\left\{ \begin{array}{ll}
\|f\|_{0,T}& \text{ if } f|_T\neq0,\\
1 & \text{ if } f|_T=0.
\end{array}
\right.
 \end{equation}
\begin{remark}
It is easy to see that
\begin{align}\label{formula16}
0\leq\xi_{T}(w_{h})\leq h_{T}, \quad \forall w_{h} \in V_{h}^{D}, \ T\in\mathcal{T}_{h}.
\end{align}
In fact, from \eqref{Ruh} and \eqref{choice-A} we have
\begin{equation}\label{formula16+}
\|R(w_{h})\|_{0,T}\leq A(w_{h})+\|f\|_{0,T},
\end{equation}
which plus the definition (\ref{formula6}) of $\xi_{T}(\cdot)$ yields (\ref{formula16}).
\end{remark}

\begin{remark}
We note that the modified version \eqref{formula6} of the artificial diffusion $\xi_{T}(\cdot)$
unifies the two cases of $ \sigma>0$ and $ \sigma=0$ into one form and allows for finite upper and lower bounds for the resulting operator $D_{T}(\cdot;\cdot,\cdot)$.
\end{remark}

\section{Existence and uniqueness  of approximation solution}

In this section, we shall show the existence and uniqueness of the discrete solution for the modified DD method.

\subsection{Preliminary results}
We give some preliminary results   for the forthcoming numerical analysis of the modified DD scheme (\ref{formulatwo}).

First, for any $v_{hb}=v_{h}+v_{b}\in V_{hb}^{D}$ and $T\in\mathcal{T}_{h}$ we easily have the orthogonality property
\begin{equation*}
\|\nabla v_{hb}\|_{0,T}^{2}=\|\nabla v_{h}\|_{0,T}^{2}+\|\nabla v_{b}\|_{0,T}^{2},
\end{equation*}
which yields
\begin{equation}\label{formula23}
|v_{h}|_{1,T}\leq|v_{hb}|_{1,T}.
\end{equation}
We also recall the result \cite[(17)]{Santos2021}
\begin{equation}\label{formula29}
\|v_{h}\|_{0,T}\leq c_{1}\|v_{hb}\|_{0,T}.
\end{equation}
As a result, by the definition \eqref{|||-T} of $|||\cdot|||_{T}$   we get
\begin{equation}
|||v_{h}|||_{T}\leq\max\{1,c_{1}\}|||v_{hb}|||_{T}.\label{formula26}
\end{equation}

The following lemma gives coercivity and boundedness results for (\ref{formulatwo}).
\begin{lemma}\label{formula17}
For any $u_{hb},v_{hb}, w_{hb}\in V_{hb}^{D}$, there hold
\begin{equation}
a(w_{hb};v_{hb},v_{hb})
\geq |||v_{hb}|||^{2} \label{2024c'}
\end{equation}
and
\begin{equation}
a(w_{hb};u_{hb},v_{hb})
\leq (M_{1}+M_2)|||u_{hb}||| \ |||v_{hb}|||,\label{formula19'}
\end{equation}
where
$$M_{2}:=\max\limits_{T\in\mathcal{T}_{h}} M_{2,T}, \quad M_{2,T}:=\left\{\begin{array}{ll}
  \min\Big\{\frac{h_{T}}{\varepsilon},\frac{c_{2}}{h_{T}\gamma}\Big\}, & \text{ if } \gamma>0,\\
  \frac{h_{T}}{\varepsilon},  & \text{ if } \gamma=0.
  \end{array}
 \right.$$
and $ M_{1}$ is the same as in \eqref{formula5}.
\end{lemma}
\begin{proof} By  the definition \eqref{formulatwo-} of $D_{T}(\cdot;\cdot,\cdot)$ and the estimate \eqref{formula16}
we immediately get, for any $w_{hb},u_{hb}, v_{hb}\in V_{hb}^{D}$,
\begin{align}
\label{formula18}
&D_{T}(w_{hb};v_{hb},v_{hb})\geq 0,
\\
&D_{T}(w_{hb};u_{hb},v_{hb})\leq \xi_{T}(w_{h})\|\nabla u_{hb}\|_{0,T}\|\nabla v_{hb}\|_{0,T}
\leq {h_{T}} |u_{hb}|_{1,T} |v_{hb}|_{1,T}. \label{formula181}
\end{align}
Then the desired coercivity result  \eqref{2024c'} follow from the coercivity \eqref{formula4} of $B(\cdot,\cdot)$.

From the Cauchy inequality and the definition \eqref{|||-T} of $|||\cdot|||_{T}$   we easily get
$$h_T|u_{hb}|_{1,T} |v_{hb}|_{1,T}\leq\frac{h_{T}}{\varepsilon}\sqrt{\varepsilon}|u_{hb}|_{1,T}\sqrt{\varepsilon}|v_{hb}|_{1,T}\leq \frac{h_{T}}{\varepsilon}|||u_{hb}|||_{T}|||v_{hb}|||_{T}.$$
When $\gamma>0$,  by the inverse inequality we also have
$$h_T|u_{hb}|_{1,T} |v_{hb}|_{1,T}\leq
\frac{c_{2}}{h_{T}\gamma}\sqrt{\gamma}\|u_{hb}\|_{0,T}\sqrt{\gamma}\|v_{hb}\|_{0,T}\leq \frac{c_{2}}{h_{T}\gamma}|||u_{hb}|||_{T}|||v_{hb}|||_{T}.$$
As a result, the above two inequalities plus \eqref{formula181} and the the boundedness \eqref{formula5} of $B(\cdot,\cdot)$  imply the desired result \eqref{formula19'}.
\end{proof}

\begin{lemma}\label{formula22}
For any $u_{hb}=u_{h}+u_{b},v_{hb}=v_{h}+v_{b}\in V_{hb}^{D}$ and $T\in\mathcal{T}_{h}$, there hold
\begin{align}
&\big| A_T(u_{h})-A_T(v_{h})\big|
\leq\sqrt{2}\max\{1,c_{1}\}M_{3,T}|||u_{hb}-v_{hb}|||_{T} ,\label{formula39'}\\
&\big|\|R(u_{h})\|_{0,T}-\|R(v_{h})\|_{0,T}\big|\leq\sqrt{2}\max\{1,c_{1}\}M_{3,T}|||u_{hb}-v_{hb}|||_{T}\label{formula24}
\end{align}
and
\begin{align}
&|\xi_{T}(u_{h})-\xi_{T}(v_{h})|\leq h_{T}M_{4,T}|||u_{hb}-v_{hb}|||_{T},\label{formula16'}
\end{align}
 where
 \begin{align}\label{M3T}
 M_{3,T}:=\left\{\begin{array}{ll}
 \max\Big\{\frac{\|\beta\|_{0,\infty,T}}{\sqrt{\varepsilon}},\frac{\|\sigma\|_{0,\infty,T}}{\sqrt{\gamma}}\Big\}, &\text{if}\ \gamma>0,\\
  \|\beta\|_{0,\infty,T}\varepsilon^{-1/2}, &\text{if}\ \gamma=0,
  \end{array}
 \right.
 \end{align}
and
 \begin{eqnarray*}
M_{4,T} :=2\sqrt{2}\max\{1,c_{1}\}M_{3,T}/\tau.
\end{eqnarray*}
\end{lemma}
\begin{proof} First,
from the triangle inequality, the Cauchy-Schwartz inequality and \eqref{|||-T} we have
\begin{eqnarray*}
&&\big|A_T(u_{h})-A_T(v_{h})\big|\\
&=&\big|\|\beta\|_{0,\infty,T}\left(\|\nabla u_{h}\|_{0,T}-\|\nabla v_{h}\|_{0,T}\right)+\|\sigma\|_{0,\infty,T}\left(\|u_h\|_{0,T}-\|v_h\|_{0,T}\right)\big|\nonumber\\
&\leq&\|\beta\|_{0,\infty,T}\|\nabla(u_{h}-v_{h})\|_{0,T}+\|\sigma\|_{0,\infty,T}\|u_h-v_h\|_{0,T}\nonumber\\
&\leq& \sqrt{2}M_{3,T}|||u_{h}-v_{h}|||_{T},\nonumber
\end{eqnarray*}
which plus \eqref{formula26} gives (\ref{formula39'}).

We next prove (\ref{formula24}).
On one hand, according to the definition of the residual $R(u_{h})$ on each element $T$, we have
\begin{eqnarray}
&\ &\big|\|R(u_{h})\|_{0,T}-\|R(v_{h})\|_{0,T}\big|_{0,T}\nonumber\\
&\leq&\|R(u_{h})-R(v_{h})\|_{0,T}\nonumber\\
&=&\|\beta\cdot\nabla(u_{h}-v_{h})+\sigma(u_{h}-v_{h})\|_{0,T}\nonumber\\
&\leq&\left(\|\beta\|_{0,\infty,T}|u_{h}-v_{h}|_{1,T}+\|\sigma\|_{0,\infty,T}\|u_{h}-v_{h}\|_{0,T}\right).
\label{formula27}
\end{eqnarray}
When $\gamma>0$ we further get
\begin{eqnarray}
&\ &\big|\|R(u_{h})\|_{0,T}-\|R(v_{h})\|_{0,T}\big|\nonumber\\
\label{formula27*}
&\leq&\sqrt{2}\max\big\{\|\beta\|_{0,\infty,T}\varepsilon^{-1/2},\|\sigma\|_{0,\infty,T}\gamma^{-1/2}\big\}|||u_{h}-v_{h}|||_{T}.
\end{eqnarray}
When $\gamma=0$, the assumption (\ref{formula2}) means $\sigma=0$, then (\ref{formula27}) leads to
\begin{eqnarray}
\big|\|R(u_{h})\|_{0,T}-\|R(v_{h})\|_{0,T}\big|\leq\sqrt{2}\|\beta\|_{0,\infty,T}\varepsilon^{-1/2}|||u_{h}-v_{h}|||_{T}.\label{december4}
\end{eqnarray}
Combining (\ref{formula27*})-(\ref{december4}) and \eqref{formula26} yields the desired estimate (\ref{formula24}).

Finally, let us show \eqref{formula16'}.

By the triangle inequality, \eqref{formula16+}, \eqref{formula39'} and \eqref{formula24} we obtain
\begin{align}
&|\xi_{T}(u_{h})-\xi_{T}(v_h)|=h_{T}\Big|\frac{\|R(u_{h})\|_{0,T}}{A_T(u_h)+\tau}-\frac{\|R(v_h)\|_{0,T}}{A_T(v_h)+\tau}\Big|\nonumber\\
\leq&\frac{h_{T}}{A_T(u_h)+\tau}\big|\|R(u_h)\|_{0,T}-\|R(v_h)\|_{0,T}\big|\nonumber\\
&\ +\frac{h_{T}\|R(v_h)\|_{0,T}}{(A_T(u_h)+\tau)(A(v_h)+\tau)}\Big|A_T(u_h)-A_T(v_h)\Big|\nonumber\\
\leq&\frac{\max\{1,c_{1}\}h_{T}}{A_T(u_h)+\tau}\big(\sqrt{2}M_{3,T}+\sqrt{2}M_{3,T}\big)|||u_{hb}-v_{hb}|||_{T}\nonumber\\
\leq& \label{formula36}
h_{T}\frac{\max\{1,c_{1}\}}{\tau}2\sqrt{2}M_{3,T}|||u_{hb}-v_{hb}|||_{T},
\end{align}
i.e. \eqref{formula16'} holds. This finishes the proof.
\end{proof}

We end this section by introducing the following  extension version of the Brouwer fixed point theorem (cf. \cite[Lemma 1.4]{Temam}):

\begin{lemma}\label{formula30}
Let $X_{h}$ be a finite dimensional Hilbert space with inner product $[\cdot,\cdot]$ and norm $[\cdot]$, and let $P : X_{h}\rightarrow X_{h}$ be a continuous operator satisfying
\begin{equation*}
[Pv,v]>0\ {\rm for}\ v\in X_{h}\ {\rm with}\ [v]=l\ {\rm for\ a\ given}\ l\in R^{+}.
\end{equation*}
Then there exists $v\in X_{h}$, with $[v]\leq l$, such that $P(v)=0$.
\end{lemma}

\subsection{Existence and uniqueness of approximation solution}

We follow  a similar routine as in \cite{Santos2021} to prove the following existence theorem:

\begin{theorem}\label{formula31}
Under the data assumption (\ref{formula2}), there exists at least one solution $u_{hb}\in V_{hb}^{D}$ to the modified DD scheme (\ref{formulatwo}) and there holds
\begin{equation} \label{bound-uh}
|||u_{hb}|||\leq
\frac{\|f\|_{0}}{\max\{\sqrt{\gamma}, \sqrt{\varepsilon}/c_0\}}+\frac{c_{3}}{\sqrt{\varepsilon}}\|g\|_{H^{-1/2}(\Gamma_{N})}=: D(f,g).
\end{equation}
\end{theorem}
\begin{proof}
Define the inner product $[\cdot,\cdot]$ on $V_{hb}^{D}\subset H_{D}^{1}(\Omega)$ by
\begin{equation*}
[u_{hb},v_{hb}]:=\varepsilon(\nabla u_{hb},\nabla v_{hb})+\gamma(u_{hb},v_{hb})\ \ \forall u_{hb}, v_{hb}\in V_{hb}^{D},
\end{equation*}
which induces a norm $[v_{hb}]=|||v_{hb}|||$. Note that this inner product is reduced to $[u_{hb},v_{hb}]=\varepsilon(\nabla u_{hb},\nabla v_{hb})$ when $\gamma=0$.
In what follows we only prove the case of $\gamma>0$, since  the case  of $\gamma=0$  is  more simpler  and the conclusion follows similarly.

Denote by $(V_{hb}^{D})^{*}$   the dual space of $V_{hb}^{D}$ with dual product $<\cdot,\cdot>$ and  by $\mathbf R : (V_{hb}^{D})^{*}\rightarrow V_{hb}^{D}$  the Riesz operator.  Define the mappings $G : V_{hb}^{D}\rightarrow (V_{hb}^{D})^{*}$ and $P_{hb} = \mathbf R\circ G: V_{hb}^{D}\rightarrow V_{hb}^{D}$    by
\begin{eqnarray}
&& [P_{hb}w_{hb},v_{hb}]=<G(w_{hb}),v_{hb}>\nonumber\\
&=& a(w_{hb};w_{hb},v_{hb})-(f,v_{hb})-(g,v_{hb})_{\Gamma_{N}} \nonumber\\
&=&B( w_{hb},v_{hb})+\sum\limits_{T\in\mathcal{T}_{h}}D_{T}(w_{hb};w_{hb},v_{hb} )-(f,v_{hb})-(g,v_{hb})_{\Gamma_{N}}, \label{formula32}
\end{eqnarray}
for any $w_{hb}, v_{hb}\in V_{hb}^{D}.$

In what follows let us show that the nonlinear operator $P_{hb}$ is continuous. In fact, from \eqref{formula32} we get, for $ u_{hb}=u_h+u_b, w_{hb}=w_h+w_b\ \in V_{hb}^{D},$
\begin{eqnarray*}
&\ &|||P_{hb}u_{hb}-P_{hb}w_{hb}|||^2=[P_{hb}u_{hb}-P_{hb}w_{hb},P_{hb}u_{hb}-P_{hb}w_{hb}]\nonumber\\
&=&a(u_{hb};u_{hb},P_{hb}u_{hb}-P_{hb}w_{hb})-a(w_{hb};w_{hb},P_{hb}u_{hb}-P_{hb}w_{hb})\nonumber\\
&=&a(u_{hb};u_{hb},P_{hb}u_{hb}-P_{hb}w_{hb})-a(w_{hb};u_{hb},P_{hb}u_{hb}-P_{hb}w_{hb})\nonumber\\
&& +a(w_{hb};u_{hb},P_{hb}u_{hb}-P_{hb}w_{hb})-a(w_{hb};w_{hb},P_{hb}u_{hb}-P_{hb}w_{hb})\nonumber\\
&=&a(w_{hb};u_{hb}-w_{hb},P_{hb}u_{hb}-P_{hb}w_{hb})\nonumber\\
&& +\sum\limits_{T\in\mathcal{T}_{h}}\int_{T}\left(\xi_{T}(u_{h})-\xi_{T}( w_{h})\right)\nabla u_{hb}\cdot\nabla(P_{hb}u_{hb}-P_{hb}w_{hb})d \mathbf x,
 \end{eqnarray*}
which, together with \eqref{formula19'} and \eqref{formula16'}, further gives
\begin{eqnarray*}
&\ &|||P_{hb}u_{hb}-P_{hb}w_{hb}|||^2 \\
&\leq&\displaystyle (M_{1}+M_2)|||u_{hb}-w_{hb}|||\ |||P_{hb}u_{hb}-P_{hb}w_{hb}|||\nonumber\\
&&\ \displaystyle+\sum\limits_{T\in\mathcal{T}_{h}}|\xi_{T}(u_{h})-\xi_{T}( w_{h}|\|\nabla u_{hb}\|_{0,T}\|\nabla(P_{hb}u_{hb}-P_{hb}w_{hb})\|_{0,T}\nonumber\\
&\leq&\displaystyle (M_{1}+M_2)|||u_{hb}-w_{hb}|||\ |||P_{hb}u_{hb}-P_{hb}w_{hb}|||\nonumber\\
&&\ \displaystyle+\sum\limits_{T\in\mathcal{T}_{h}}h_{T}M_{4,T}|||u_{hb}-w_{hb}|||_{T}\|\nabla u_{hb}\|_{0,T}\|\nabla(P_{hb}u_{hb}-P_{hb}w_{hb})\|_{0,T}\nonumber\\
&\leq&\displaystyle (M_{1}+M_2)|||u_{hb}-w_{hb}|||\ |||P_{hb}u_{hb}-P_{hb}w_{hb}|||\nonumber\\
&&\ \displaystyle+h\frac{M_{4}}{\sqrt{\varepsilon}}\sum\limits_{T\in\mathcal{T}_{h}}|||u_{hb}-w_{hb}|||_{T}\|\nabla u_{hb}\|_{0,T}|||P_{hb}u_{hb}-P_{hb}w_{hb}|||_{T}\nonumber\\
&\leq&\displaystyle\left(M_{1}+M_{2}+h\frac{\|\nabla u_{hb}\|_{0}}{\sqrt{\varepsilon}}M_{4}\right)|||u_{hb}-w_{hb}||| \  |||P_{hb}u_{hb}-P_{hb}w_{hb}|||,\nonumber
 \end{eqnarray*}
where
\begin{equation*}
M_{4}=\max\limits_{T\in\mathcal{T}_{h}}M_{4,T}.
\end{equation*}
Note that we utilize the Cauchy inequality twice in the last step above.

This means that
\begin{eqnarray*}
&\ &|||P_{hb}u_{hb}-P_{hb}w_{hb}|||
\leq  \left(M_{1}+M_{2}+h\frac{\|\nabla u_{hb}\|_{0}}{\sqrt{\varepsilon}}M_{4}\right)|||u_{hb}-w_{hb}|||.
\end{eqnarray*}
Notice that
\begin{equation*}
M_{2,T}\leq h_{T}/\varepsilon,
\end{equation*}
which means $M_{2}\leq h/\varepsilon$. As a consequence,
\begin{equation*}
M_{1}+M_{2}+h\frac{\|\nabla u_{hb}\|_{0}}{\sqrt{\varepsilon}}M_{4}
\end{equation*}
is bounded even if the mesh size $h$ approaches to zero, and the continuity of the operator $P_{hb}$ follows.

Finally, let us prove that
$[P_{hb}u_{hb},u_{hb}]>0$ if $[u_{hb}]=|||u_{hb}|||=l$  for a sufficiently large $l$. Notice that
\begin{equation}\label{formula42}
[P_{hb}u_{hb},u_{hb}]=a(u_{hb};u_{hb},u_{hb})-(f,u_{hb})-(g,u_{hb})_{\Gamma_{N}}.
\end{equation}
In light of $u_{hb}|_{\Gamma_{D}}=0$  and the Poincar\'{e} inequality  we have
\begin{eqnarray*}
|(g,u_{hb})_{\Gamma_{N}}|&\leq&\|g\|_{H^{-1/2}(\Gamma_{N})}\|u_{hb}\|_{H^{1/2}(\Gamma_{N})}\nonumber\\
&\leq&\|g\|_{H^{-1/2}(\Gamma_{N})}\|u_{hb}\|_{1}\nonumber\\
&\leq&c_{3}\|g\|_{H^{-1/2}(\Gamma_{N})}|u_{hb}|_{1}\nonumber\\
&\leq&\frac{c_{3}}{\sqrt{\varepsilon}}\|g\|_{H^{-1/2}(\Gamma_{N})}|||u_{hb}|||,\label{duJ1}
\end{eqnarray*}
which, together with (\ref{formula42}) and \eqref{2024c'}, leads to
\begin{eqnarray*}
[P_{hb}u_{hb},u_{hb}]&\geq&|||u_{hb}|||^{2}-(f,u_{hb})-(g,u_{hb})_{\Gamma_{N}}\nonumber\\
&\geq&|||u_{hb}|||^{2}-\frac{1}{\sqrt{\gamma}}\|f\|_{0}|||u_{hb}|||-\frac{c_{3}}{\sqrt{\varepsilon}}\|g\|_{H^{-1/2}(\Gamma_{N})}|||u_{hb}|||\nonumber\\
&=&l\left(l-\frac{1}{\sqrt{\gamma}}\|f\|_{0}-\frac{c_{3}}{\sqrt{\varepsilon}}\|g\|_{H^{-1/2}(\Gamma_{N})}\right)\\
&>&0,
\end{eqnarray*}
provided that $l>\frac{1}{\sqrt{\gamma}}\|f\|_{0}+\frac{c_{3}}{\sqrt{\varepsilon}}\|g\|_{H^{-1/2}(\Gamma_{N})}$.

As a result, the operator $P_{hb}$ satisfies the conditions of Lemma \ref{formula30}, then there exists $u_{hb}\in V_{hb}^{D}$, such that $P_{hb}u_{hb}=0$, which plus
 (\ref{formula32}) yields
\begin{equation}\label{duJ2}
\displaystyle a(u_{hb};u_{hb},v_{hb})=(f,v_{hb})+(g,v_{hb})_{\Gamma_{N}}, \forall\ v_{hb}\in V_{hb}^{D}.
\end{equation}
This means that the modified DD scheme (\ref{formulatwo}) admits a solution.

Taking $v_{hb}=u_{hb}$ in \eqref{duJ2}, we have from \eqref{2024c'} and \eqref{duJ1}
\begin{eqnarray}
|||u_{hb}|||^{2}&\leq&a(u_{hb};u_{hb},u_{hb})\nonumber\\
&\leq&\|f\|_{0}\|u_{hb}\|_{0}+\|g\|_{H^{-1/2}(\Gamma_{N})}\|u_{hb}\|_{H^{1/2}(\Gamma_{N})}\nonumber\\
&\leq&\frac{\|f\|_{0}}{\max\{\sqrt{\gamma},\sqrt{\varepsilon}/c_0\}}|||u_{hb}|||+\frac{c_{3}}{\sqrt{\varepsilon}}\|g\|_{H^{-1/2}(\Gamma_{N})}|||u_{hb}|||\nonumber\\
&\leq&\left(\frac{\|f\|_{0}}{\max\{\sqrt{\gamma},\sqrt{\varepsilon}/c_0\}}+\frac{c_{3}}{\sqrt{\varepsilon}}\|g\|_{H^{-1/2}(\Gamma_{N})}\right)|||u_{hb}|||,\nonumber
\end{eqnarray}
which results in the desired estimate \eqref{bound-uh}. This completes the proof.
\end{proof}

\begin{remark}
In the proof of Theorem 3.1 in \cite{Santos2021}, the boundedness of $\|\nabla u_{hb}\|_{0,\infty}$ ($\forall u_{hb}\in V_{hb}^{D}$) has been used. According to the inverse
inequality, it holds
$$\|\nabla u_{hb}\|_{0,\infty}\leq c_2h^{-d/2}\|\nabla u_{hb}\|_{0}.$$
As a result, $\|\nabla u_{hb}\|_{0,\infty}$ will approach to infinity when the mesh size $h$ tends to zero. However, {
we only assume the boundedness of $\|\nabla u_{hb}\|_{0}$ in our analysis}.
\end{remark}

\begin{theorem}\label{uniqness-uh}
Under the same conditions as in Theorem \ref{formula31} and the smallness condition
\begin{equation} \label{smallness-uh}
h<\varepsilon\left(\frac{\|f\|_{0}}{\max\{\sqrt{\gamma}, \sqrt{\varepsilon}/c_0\}}+\frac{c_{3}}{\sqrt{\varepsilon}}\|g\|_{H^{-1/2}(\Gamma_{N})}\right)^{-1}M_{4}^{-1},
\end{equation}
the modified DD scheme (\ref{formulatwo}) admits a unique solution $u_{hb}\in V_{hb}^{D}$.
\end{theorem}
\begin{proof}
Assume that there are two solutions, $u_{hb}^1, u_{hb}^2\in V_{hb}^{D}$, to (\ref{formulatwo}), i.e. for $i=1,2$,
\begin{equation*}
\displaystyle a(u_{hb}^i;u_{hb}^i,v_{hb})=(f,v_{hb})+(g,v_{hb})_{\Gamma_{N}}, \forall\ v_{hb}\in V_{hb}^{D}.
\end{equation*}
Then we have
\begin{eqnarray*}
0&=&a(u_{hb}^1;u_{hb}^1,u_{hb}^1-u_{hb}^2)-a(u_{hb}^2;u_{hb}^2,u_{hb}^1-u_{hb}^2)\\
&=&a(u_{hb}^1;u_{hb}^1-u_{hb}^2,u_{hb}^1-u_{hb}^2)\nonumber\\
&& +\sum\limits_{T\in\mathcal{T}_{h}}\int_{T}\left(\xi_{T}(u_{h}^1)-\xi_{T}( u_{h}^2)\right)\nabla u_{hb}^2\cdot\nabla (u_{hb}^1-u_{hb}^2)d \mathbf x,
\end{eqnarray*}
which, together with \eqref{2024c'} and \eqref{formula16'}, gives
\begin{eqnarray}
|||u_{hb}^1-u_{hb}^2|||^2&\leq&  a(u_{hb}^1;u_{hb}^1-u_{hb}^2,u_{hb}^1-u_{hb}^2)\nonumber\\
&\leq&\Big|\sum\limits_{T\in\mathcal{T}_{h}}\int_{T}\left(\xi_{T}(u_{h}^1)-\xi_{T}( u_{h}^2)\right)\nabla u_{hb}^2\cdot\nabla (u_{hb}^1-u_{hb}^2)d\mathbf x\Big|\nonumber\\
&\leq&\sum\limits_{T\in\mathcal{T}_{h}}h_{T}M_{4,T}|||u_{hb}^1-u_{hb}^2|||_{T}\|\nabla u_{hb}^2\|_{0,T}\|\nabla(u_{hb}^1-u_{hb}^2)\|_{0,T}\nonumber\\
&\leq&h\varepsilon^{-1}M_{4}\sum\limits_{T\in\mathcal{T}_{h}}|||u_{hb}^1-u_{hb}^2|||_{T}|||u_{hb}^2|||_{T}|||u_{hb}^1-u_{hb}^2|||_{T}\nonumber\\
&\leq&h\varepsilon^{-1}M_{4}|||u_{hb}^2||| \ |||u_{hb}^1-u_{hb}^2|||^{2}.\label{duJ3}
\end{eqnarray}

Applying the stability estimation \eqref{bound-uh} to \eqref{duJ3}, we have
\begin{equation*}
|||u_{hb}^1-u_{hb}^2|||^2\leq\frac{hM_{4}}{\varepsilon}\left(\frac{\|f\|_{0}}{\max\{\sqrt{\gamma},\sqrt{\varepsilon}/c_0\}}+
\frac{c_{3}}{\sqrt{\varepsilon}}\|g\|_{H^{-1/2}(\Gamma_{N})}\right)|||u_{hb}^1-u_{hb}^2|||^{2}.
\end{equation*}
In view of the smallness condition \eqref{smallness-uh} and the above inequality, we know that $u_{hb}^1=u_{hb}^2.$
\end{proof}

\section{Optimal estimate}

For the DD method with the artificial diffusion (\ref{formula6-}) or (\ref{formula7-}),
a suboptimal convergence rate   $O(h^{1/2})$ has been derived in energy norm in the literature \cite{Santos2021}, though  the optimal convergence rate $O(h)$ being observed in  numerical experiments. In this section, we shall establish  the first-order optimal convergence rate for the modified DD method.

\begin{theorem}\label{formula45}
Under the assumption (\ref{formula2}), let $u\in H^{2}(\Omega)\cap H_{D}^{1}(\Omega)$ be the solution of the continuous problem \eqref{formula1} and let $u_{hb}\in V_{hb}^{D}$ be the solution of the modified DD scheme (\ref{formulatwo}). 
Then there holds
\begin{align}
 &\displaystyle |||u-u_{hb}|||^{2}+\sum\limits_{T\in\mathcal{T}_{h}}D_{T}(u_{hb};u_{hb},u_{hb})\nonumber\\
\leq &
c_{5}^{2}\left(4\gamma h^{2}+4 \varepsilon+ 2M_{5} + 2h\right)h^{2}\|u\|_{2}^{2}  +2\left(\|f\|_{0}+c_{6}M_{3}D(f,g)\right)\frac{\sqrt{d}}{\varepsilon}h^{d/2}h^{2}|u|_{1}, \label{formula46}
\end{align}
which further gives
\begin{align}
 \displaystyle |||u-u_{hb}|||
 \leq &\ \
c_{5}\left(4\gamma h^{2}+4 \varepsilon+ 2M_{5} + 2h\right)^\frac12 h \|u\|_{2} \nonumber \\
&\quad   +\left(2\sqrt{d}\|f\|_{0}+2c_{6}\sqrt{d}M_{3}D(f,g)\right)^\frac12 \frac{h^{d/4}}{\sqrt{\varepsilon}  }h|u|_{1}^\frac12, \label{formula466}
\end{align}
where
\begin{equation*}
M_{5}:=\min\left\{\frac{1 }{ \varepsilon}\left(\varepsilon    +c_0\|\beta\|_{0,\infty}+c_0\|\sigma\|_{0,\infty} h \right)^2,   \varepsilon+\frac{ 1}{\gamma}(\|\beta\|_{0,\infty}
+  \|\sigma\|_{0,\infty}h )^2 \right \}.
\end{equation*}
\end{theorem}
\begin{proof}
Denote by $\pi_{h} : H_{D}^{1}(\Omega)\cap H^{2}(\Omega)\rightarrow V_{h}^{D}$ the  usual Lagrange interpolation operator, and recall its standard interpolation estimate (\cite{Ciarlet})
\begin{equation}\label{formula47}
\|v-\pi_{h}v\|_{0,T}+h_{T}\|\nabla(v-\pi_{h}v)\|_{0,T}\leq c_{5}h_{T}^{2}\|v\|_{2,T},\ \ \forall v\in H_{D}^{1}(\Omega)\cap H^{2}(\Omega)
\end{equation}
for any $T\in\mathcal{T}_{h}$,
where $c_{5}$ is a positive constant depending on the shape regularity of $T$.

Set
\begin{equation}\label{formula48}
u-u_{hb}=\eta+\varphi \quad \text{with } \eta:=u-\pi_{h}u,\varphi:=\pi_{h}u-u_{hb}.
\end{equation}
Applying (\ref{formula47}) gives
\begin{equation}\label{formula49}
|||\eta|||=\left(\varepsilon|u-\pi_{h}u|_{1}^{2}+\gamma\|u-\pi_{h}u\|_{0}^{2}\right)^{1/2}\leq c_{5}(\gamma^{1/2} h^{2}+\varepsilon^{1/2}h)\|u\|_{2}.
\end{equation}

Noticing $\varphi\in V_{hb}^{D}\subset H_{D}^{1}(\Omega)$, from (\ref{formula3}) and (\ref{formulatwo}) we have
\begin{equation*}\label{formula50}
\displaystyle B(u-u_{hb},\varphi)=\sum\limits_{T\in\mathcal{T}_{h}}D_{T}(u_{hb};u_{hb},\varphi),
\end{equation*}
which plus the coercivity (\ref{formula4}) yields
\begin{align} \label{formula51}
&\displaystyle|||\varphi|||^{2}+\sum\limits_{T\in\mathcal{T}_{h}}D_{T}(u_{hb};u_{hb},u_{hb})\nonumber\\
\leq&\displaystyle B(\varphi,\varphi)+\sum\limits_{T\in\mathcal{T}_{h}}D_{T}(u_{hb};u_{hb},u_{hb})\nonumber\\
=&B(\pi_{h}u-u,\varphi)+B(u-u_{hb},\varphi)+\sum\limits_{T\in\mathcal{T}_{h}}D_{T}(u_{hb};u_{hb},\pi_{h}u-\varphi)\nonumber\\
=&\displaystyle-B(\eta,\varphi)+\sum\limits_{T\in\mathcal{T}_{h}}D_{T}(u_{hb};u_{hb},\pi_{h}u).
\end{align}
For the  term $-B(\eta,\varphi)$, we easily have
\begin{align}
 |-B(\eta,\varphi)|
\leq\varepsilon\|\nabla\eta\|_{0}\|\nabla\varphi\|_{0}+\|\beta\|_{0,\infty}\|\nabla\eta\|_{0}\|\varphi\|_{0}
+\|\sigma\|_{0,\infty}\|\eta\|_{0}\|\varphi\|_{0}.
\end{align}
Noticing that
$|||\varphi|||^{2} =\varepsilon|\varphi|_{1}^{2}+\gamma\|\varphi\|_{0}^{2}$ by \eqref{|||-norm}, we apply the Yang's inequality, the inequality \eqref{formula0-1} and the estimate (\ref{formula47}) to obtain  
\begin{eqnarray*}
  |-B(\eta,\varphi)|
& \leq&\left(\varepsilon\|\nabla\eta\|_{0} +c_0\|\beta\|_{0,\infty}\|\nabla\eta\|_{0}
+c_0\|\sigma\|_{0,\infty}\|\eta\|_{0})\right\|\nabla\varphi\|_{0}\nonumber\\
& \leq &\frac{1}{2\varepsilon}\left(\varepsilon\|\nabla\eta\|_{0} +c_0\|\beta\|_{0,\infty}\|\nabla\eta\|_{0}
+c_0\|\sigma\|_{0,\infty}\|\eta\|_{0}\right)^2+\frac\varepsilon 2 \|\nabla\varphi\|_{0}^2\nonumber\\
&\leq &   \frac{c_5^2 }{ 2\varepsilon}\left(\varepsilon    +c_0\|\beta\|_{0,\infty}+c_0\|\sigma\|_{0,\infty} h \right)^2
 h^2\|u\|_{2}^2+\frac{1}{2}|||\varphi|||^{2}.
\end{eqnarray*}
In particular, when $\gamma >0$ we also have
\begin{align*}
  |-B(\eta,\varphi)|&\leq\varepsilon\|\nabla\eta\|_{0}\|\nabla\varphi\|_{0}+\left(\|\beta\|_{0,\infty}\|\nabla\eta\|_{0}
+\|\sigma\|_{0,\infty}\|\eta\|_{0}\right)\|\varphi\|_{0}\nonumber\\
& \leq \frac{\varepsilon}{2}\|\nabla\eta\|_{0}^{2}+\frac{\varepsilon}{2}\|\nabla\varphi\|_{0}^{2}+\frac{1}{2\gamma}\left(\|\beta\|_{0,\infty}\|\nabla\eta\|_{0}
+\|\sigma\|_{0,\infty}\|\eta\|_{0}\right)^2+\frac{\gamma}{2}\|\varphi\|_{0}^{2}\nonumber\\
& \leq \frac12c_5^2h^2\left( \varepsilon+\frac{ 1}{\gamma}(\|\beta\|_{0,\infty}
+  \|\sigma\|_{0,\infty}h )^2\right)\|u\|_{2}^2 +\frac{1}{2}|||\varphi|||^{2}. 
\end{align*}
As a result, we have
\begin{eqnarray}\label{formula52}
  |-B(\eta,\varphi)| 
 \leq \frac{1}{2}c_5^2M_5 h^{2}\|u\|_{2}^{2}+\frac{1}{2}|||\varphi|||^{2},
\end{eqnarray}
where
\begin{equation*}
M_{5}=\min\left\{\frac{1 }{ \varepsilon}\left(\varepsilon    +c_0\|\beta\|_{0,\infty}+c_0\|\sigma\|_{0,\infty} h \right)^2,   \varepsilon+\frac{ 1}{\gamma}(\|\beta\|_{0,\infty}
+  \|\sigma\|_{0,\infty}h )^2 \right \}.
\end{equation*}

For the term $\displaystyle\sum\limits_{T\in\mathcal{T}_{h}}D_{T}(u_{hb};u_{hb},\pi_{h}u)$ on the right side of \eqref{formula51}, noticing that $$u_{hb}=u_{h}+u_{b}, \quad \kappa_{h}(u_{hb})=u_{h}$$ and
\begin{equation*}
\displaystyle\int_{T}\nabla v_{h}\cdot\nabla v_{b}d{\bf x}=0,\ \forall\ v_{h}\in V_{h},\ \forall\ v_{b}\in V_{b},\forall\ T\in\mathcal{T}_{h},
\end{equation*}
we obtain
\begin{eqnarray}
&\ &\displaystyle\sum\limits_{T\in\mathcal{T}_{h}}D_{T}(u_{hb};u_{hb},\pi_{h}u)\nonumber\\
&=&\sum\limits_{T\in\mathcal{T}_{h}}\int_{T}\xi_{T}(\kappa_{h}(u_{hb}))\nabla u_{hb}\cdot\nabla(\pi_{h}u)d{\bf x}\nonumber\\
&=&\sum\limits_{T\in\mathcal{T}_{h}}\xi_{T}(u_{h})\int_{T}\nabla u_{h}\cdot\nabla(\pi_{h}u)d{\bf x}\nonumber\\
&=&\sum\limits_{T\in\mathcal{T}_{h}}\xi_{T}(u_{h})\int_{T}\nabla u_{h}\cdot\nabla(\pi_{h}u-u)
+\sum\limits_{T\in\mathcal{T}_{h}}\xi_{T}(u_{h})\int_{T}\nabla u_{h}\cdot\nabla u.\label{formula53}
\end{eqnarray}

For the first term on the right hand side of the last line of (\ref{formula53}), in light of (\ref{formula47}), (\ref{formula6}) and (\ref{formula23})
we get
\begin{eqnarray}
&  &  \Big|\sum\limits_{T\in\mathcal{T}_{h}}\xi_{T}(u_{h})\int_{T}\nabla u_{h}\cdot\nabla(\pi_{h}u-u)d{\bf x}\Big|\nonumber\\
&\leq&\sum\limits_{T\in\mathcal{T}_{h}}\xi_{T}(u_{h})\|\nabla u_{h}\|_{0,T}\|\nabla(u-\pi_{h}u)\|_{0,T}\nonumber\\
&\leq&c_{5}\sum\limits_{T\in\mathcal{T}_{h}}\xi_{T}(u_{h})\|\nabla u_{h}\|_{0,T}h_{T}\|u\|_{2,T}\nonumber\\
&\leq&c_{5}h\max_{T\in\mathcal{T}_{h}}\sqrt{\xi_{T}(u_{h})}\|u\|_{2}\Big\{\sum\limits_{T\in\mathcal{T}_{h}}\xi_{T}(u_{h})\|\nabla u_{hb}\|_{0,T}^{2}\Big\}^{1/2}\nonumber\\
\label{formula54}
&\leq&c_{5}h^{3/2}\Big\{\sum\limits_{T\in\mathcal{T}_{h}}\xi_{T}(u_{h})\|\nabla u_{hb}\|_{0,T}^{2}\Big\}^{1/2}\|u\|_{2},
\end{eqnarray}
where in the last ``$\leq$" we have employed the estimate \eqref{formula16}.

Next, let us  estimate the second term on the right hand side of (\ref{formula53}).
According to the choice (\ref{formula6}), the artificial diffusion coefficient $\xi_{T}(u_{h})$ vanishes when the local Pecl\'{e}t number $P_{e_{T}}\leq1$. Notice that
$P_{e_{T}}>1$ indicates $\|\beta\|_{0,\infty,T}>0$. As a result, we have
\begin{equation}\label{Julyduone}
\displaystyle\sum\limits_{T\in\mathcal{T}_{h}}\xi_{T}(u_{h})\int_{T}\nabla u_{h}\cdot\nabla u d{\bf x}=\sum\limits_{T\in\mathcal{T}_{h},P_{e_{T}}>1}\xi_{T}(u_{h})\int_{T}\nabla u_{h}\cdot\nabla u d{\bf x},
\end{equation}
which, together with 
(\ref{formula6})-\eqref{choice-tau}, \eqref{formula16+}, \eqref{M3T} 
and (\ref{formula26}),  yields
\begin{eqnarray}
&\ &\Big|\sum\limits_{T\in\mathcal{T}_{h}}\xi_{T}(u_{h})\int_{T}\nabla u_{h}\cdot\nabla u d{\bf x}\Big|\nonumber\\
&\leq&\sum\limits_{T\in\mathcal{T}_{h},P_{e_{T}}>1}\frac{h_{T}\|R(u_{h})\|_{0,T}}{A_{T}(u_{h})+\tau}\|\nabla u_{h}\|_{0,T}|u|_{1,T}\nonumber\\
&\leq&\sum\limits_{T\in\mathcal{T}_{h},P_{e_{T}}>1}h_{T}\frac{\|f\|_{0,T}+\|\beta\|_{0,\infty,T}|u_{h}|_{1,T}+
\|\sigma\|_{0,\infty,T}\|u_{h}\|_{0,T}}{\|\beta\|_{0,\infty,T}}|u|_{1,T}\nonumber\\
&\leq&\sum\limits_{T\in\mathcal{T}_{h},P_{e_{T}}>1}h_{T}\frac{\|f\|_{0,T}+\sqrt{2}M_{3,T}|||u_{h}|||_{T}}{\|\beta\|_{0,\infty,T}}|u|_{1,T}\nonumber\\
&\leq&\sum\limits_{T\in\mathcal{T}_{h},P_{e_{T}}>1}h_{T}\frac{\|f\|_{0,T}+\sqrt{2}\max\{1,c_{1}\}M_{3,T}|||u_{hb}|||_{T}}{\|\beta\|_{0,\infty,T}}|u|_{1,T}. \label{formula56}
\end{eqnarray}

The definition of the local Pecl\'{e}t number $P_{e_{T}}$  yields
\begin{equation*}
1<P_{e_{T}}=\frac{h_{T}\|\beta\|_{0,T}}{2\varepsilon}\leq\frac{h_{T}\sqrt{d}\|\beta\|_{0,\infty,T}|T|^{1/2}}{2\varepsilon},
\end{equation*}
which implies
\begin{equation*}\label{August1}
\frac{1}{\|\beta\|_{0,\infty,T}}\leq\frac{h_{T}\sqrt{d}|T|^{1/2}}{2\varepsilon}.
\end{equation*}
Inserting this inequality into (\ref{formula56}) and applying the stability estimate (\ref{bound-uh}) of the discrete solution $u_{hb}$, we attain
\begin{eqnarray}
&\ &\Big|\sum\limits_{T\in\mathcal{T}_{h}}\xi_{T}(u_{h})\int_{T}\nabla u_{h}\cdot\nabla u d{\bf x}\Big|\nonumber\\
&\leq&\sum\limits_{T\in\mathcal{T}_{h},P_{e_{T}}>1}h_{T}\left(\|f\|_{0,T}+
\sqrt{2}\max\{1,c_{1}\}M_{3,T}|||u_{hb}|||_{T}\right)\frac{h_{T}\sqrt{d}|T|^{1/2}}{2\varepsilon}|u|_{1,T}\nonumber\\
&\leq&\sum\limits_{T\in\mathcal{T}_{h}}h_{T}^{2}\left(\|f\|_{0,T}+
\sqrt{2}\max\{1,c_{1}\}M_{3,T}|||u_{hb}|||_{T}\right)\frac{\sqrt{d}|T|^{1/2}}{2\varepsilon}|u|_{1,T}\nonumber\\
&\leq&\sum\limits_{T\in\mathcal{T}_{h}}h_{T}^{2}\left(\|f\|_{0,T}+
\sqrt{2}\max\{1,c_{1}\}M_{3,T}|||u_{hb}|||_{T}\right)\frac{\sqrt{d}h_{T}^{d/2}}{2\varepsilon}|u|_{1,T}\nonumber\\
&\leq&h^{2}\left(\|f\|_{0}+c_{6}M_{3}|||u_{hb}|||\right)\frac{\sqrt{d}}{2\varepsilon}h^{d/2}|u|_{1}\nonumber\\
&\leq&h^{2}\left(\|f\|_{0}+c_{6}M_{3}D(f,g)\right)\frac{\sqrt{d}}{2\varepsilon}h^{d/2}|u|_{1}. \label{August2}
\end{eqnarray}
where
\begin{equation*}
M_{3}:=\max_{T\in\mathcal{T}_{h}}M_{3,T}, \quad c_{6}:=\sqrt{2}\max\{1,c_{1}\},
\end{equation*}
and $D(f,g)$ is defined in \eqref{bound-uh}.

Collecting (\ref{formula54}) and (\ref{August2}), from (\ref{formula53}) we get
\begin{align}
\displaystyle &\Big|\sum\limits_{T\in\mathcal{T}_{h}}D_{T}(u_{hb};u_{hb},\pi_{h}u)\Big|\nonumber\\
\leq& c_{5}h^{\frac{3}{2}}\Big\{\sum\limits_{T\in\mathcal{T}_{h}}\xi_{T}(u_{h})\|\nabla u_{hb}\|_{0,T}^{2}\Big\}^{\frac{1}{2}}\|u\|_{2}
+h^{2}\left(\|f\|_{0}+c_{6}M_{3}D(f,g)\right)\frac{\sqrt{d}}{2\varepsilon}h^{d/2}|u|_{1}. \label{formula57}
\end{align}

In view of
\begin{equation*}
\xi_{T}(u_{h})\|\nabla u_{hb}\|_{0,T}^{2}=D_{T}(u_{hb};u_{hb},u_{hb}),
\end{equation*}
combining (\ref{formula51}), (\ref{formula52}) and (\ref{formula57}) gives
\begin{eqnarray*}
&\ &|||\varphi|||^{2}+\sum\limits_{T\in\mathcal{T}_{h}}D_{T}(u_{hb};u_{hb},u_{hb})\nonumber\\
&\leq&\frac{1}{2}c_5^2M_{5}h^{2}\|u\|_{2}^{2}+\frac{1}{2}|||\varphi|||^{2}+\frac{1}{2}\sum\limits_{T\in\mathcal{T}_{h}}D_{T}(u_{hb};u_{hb},u_{hb})\nonumber\\
&\ &\ +\frac12c_{5}^{2}h^{3}\|u\|_{2}^{2}+h^{2}\left(\|f\|_{0}+c_{6}M_{3}D(f,g)\right)\frac{\sqrt{d}}{2\varepsilon}h^{d/2}|u|_{1},\nonumber\\
\end{eqnarray*}
which yields
\begin{eqnarray}
&\ &|||\varphi|||^{2}+\sum\limits_{T\in\mathcal{T}_{h}}D_{T}(u_{hb};u_{hb},u_{hb})\nonumber\\
&\leq&\left( c_5^2M_{5}+c_{5}^{2}h\right)h^{2}\|u\|_{2}^{2}+\left(\|f\|_{0}+c_{6}M_{3}D(f,g)\right)\frac{\sqrt{d}}{\varepsilon}h^{d/2}h^{2}|u|_{1}.\nonumber
\end{eqnarray}
This estimate plus the triangle inequality, (\ref{formula48}) and (\ref{formula49})
indicates
\begin{align*}
 & |||u-u_{hb}|||^{2}+\sum\limits_{T\in\mathcal{T}_{h}}D_{T}(u_{hb};u_{hb},u_{hb})\nonumber\\
\leq&\left(|||\eta|||+|||\varphi|||\right)^2+\sum\limits_{T\in\mathcal{T}_{h}}D_{T}(u_{hb};u_{hb},u_{hb})\nonumber\\
\leq&2\left(|||\eta|||^{2}+|||\varphi|||^{2}\right)+\sum\limits_{T\in\mathcal{T}_{h}}D_{T}(u_{hb};u_{hb},u_{hb})\nonumber\\
\leq&4c_{5}^{2}(\gamma h^{4}+\varepsilon h^{2})\|u\|_{2}^{2}+2\left(|||\varphi|||^{2}+\sum\limits_{T\in\mathcal{T}_{h}}D_{T}(u_{hb};u_{hb},u_{hb})\right)\nonumber\\
\leq&c_{5}^{2}\left(4\gamma h^{2}+4 \varepsilon+2 M_{5} +2 h\right)h^{2}\|u\|_{2}^{2}  +2\left(\|f\|_{0}+c_{6}M_{3}D(f,g)\right)\frac{\sqrt{d}}{\varepsilon}h^{d/2}h^{2}|u|_{1}, 
\end{align*}
i.e. the desired estimate \eqref{formula46} holds. This completes the proof.
\end{proof}

\begin{remark}
Theorem \ref{formula45} shows that  for the modified DD method 
 the numerical error in the energy norm, i.e. $|||u-u_{hb}|||  $, is of optimal first-order convergence rate with respect to the mesh-size $h$. We point out that the key to the improvement of the theoretical accuracy in our analysis, compared with \cite{Santos2021}, lies in
that we have followed a different way to estimate the term $\displaystyle\sum\limits_{T\in\mathcal{T}_{h}}D_{T}(u_{hb};u_{hb},\pi_{h}u)$ in (\ref{formula51}).
 \end{remark}

\begin{remark}\label{pinglun2}
Notice that the upper bound of the numerical error, i.e. the right hand side of (\ref{formula466}), consists of
two terms,
$$T_1:=c_{5}\left(4\gamma h^{2}+4 \varepsilon+ 2M_{5} + 2h\right)^\frac12 h \|u\|_{2} $$
and
$$ T_2:=\left(2\sqrt{d}\|f\|_{0}+2c_{6}\sqrt{d}M_{3}D(f,g)\right)^\frac12 \frac{h^{d/4}}{\sqrt{\varepsilon}  }h|u|_{1}^\frac12. $$
The former term $T_1$ is due to the approximation error, and the latter term $T_2$ is owing to the additional artificial dissipation. To clarify the  dependence of these two terms  on $\varepsilon$ for $0<\varepsilon <<1$, we
assume that $\|f\|_{0}$ and $\|g\|_{H^{-1/2}(\Gamma_{N})}$
are independent of $\varepsilon$, $\beta \neq 0$ and   the exact solution $u$ has a layer of width $O(\varepsilon^{\alpha}) $ 
with
$$ \|u\|_{2}=\mathcal{O}(\varepsilon^{-2\alpha}), \quad |u|_{1}=\mathcal{O}(\varepsilon^{-\alpha}).$$ 
Thus, we have $M_{5}=\mathcal O(1)$
if  $\gamma>0$ and  $M_{5}=\mathcal O(\varepsilon^{-1})$ if $\gamma=0$.
As a result, there hold
$$T_1=\left\{\begin{aligned}
\mathcal O( \varepsilon^{-2\alpha}h) & \text{ if }\gamma>0, \\
\mathcal O( \varepsilon^{-2\alpha-\frac12}h)& \text{ if }\gamma=0.
\end{aligned}
\right., \quad T_2=\mathcal O( \varepsilon^{-1-\frac{\alpha}{2}}h^{ \frac d4}h).
$$
We can see that  in case of $\gamma>0$ (or $\gamma=0$), $T_{1}$ is a higher order term with respect to $ \varepsilon^{-1}$
than $T_{2}$ when $\alpha\geq2/3$ (or $\alpha\geq1/3$) due to $1+\alpha/2\leq2\alpha$ (or $1+\alpha/2\leq2\alpha+1/2$).
This suggests that the estimate component $T_{1}$ dominates under these circumstances, and that the global estimate $T_{1}+T_{2}$
is sharp with respect to $h$.

\begin{remark}\label{pinglun3}
As shown in Remark \ref{pinglun2},  the  layers may deteriorate the numerical accuracy when the mesh size $h$ is not small enough. To   resolve the boundary or interior layers it is known that   the adaptive mesh refinement    based on a posteriori error estimators is   a basic strategy  \cite{Babuska1978,DuXie,cai1,Ani,NUMER,Verfurth1996,Du2015,Du2021,Du2022}. We mention that
in a recent work \cite{Du2023} we have developed a robust a posteriori error estimator for the DD method  and  proposed a linearized adaptive dynamic diffusion algorithm to solve the related nonlinear DD system.
\end{remark}
\end{remark}

\section{Numerical experiments}
In this section, we first give a  fixed point iterative algorithm  for the  modified DD method  (\ref{formulatwo}), then perform   numerical experiments to verify the theoretical analysis.
\subsection{Fixed point iteration}

Since the scheme (\ref{formulatwo}) is a nonlinear system, we follow the   fixed point iteration process developed in \cite{Santos2018,Santos2021} to get the approximation solution:
\begin{itemize}
\item[Step 1] Use the classical SUPG method to get an initial data $u_{hb}^{0}=u_{h}^{0}$, set the initial artificial
diffusion $\xi_T^{0}(u_{h}^{0})=0$  for any $T\in\mathcal{T}_{h}$,  the maximum number, $M=30$, of nonlinear iterations and the tolerance error $\delta=10^{-6}$, and set the iteration counter $k=0$;

\item[Step 2] Apply an accelerated algorithm described in \cite{Santos2021,Santos2007} to speed up the convergence. For any $T\in\mathcal{T}_{h}$,    
compute the artificial
diffusion $\xi_T^{k+1}(u_{h}^{k+1})$ at the $k+1-$th iteration by
\begin{equation*}
\xi_T^{k+1}(u_{h}^{k+1})=\omega {\xi}_T(u_{h}^{k})+(1-\omega)\xi_T^{k}(u_{h}^{k}),
\end{equation*}
where $ {\xi}_T(\cdot)$ is given by (\ref{formula6}), and  the damping factor $\omega$ is determined by a freezing strategy, i.e.
$$\omega=\left\{ \begin{aligned}
&0  \quad  \text{ if }  {\left|\|R(u_{h}^{k-1})\|_{0,T}-\|R(u_{h}^{k})\|_{0,T}\right|}\leq 0.2 {\|R(u_{h}^{k-1})\|_{0,T}} \text{ and } k\geq 1,\\
&0.5 \quad \text{ otherwise;}
\end{aligned}
\right. $$
\item[Step 3] Solve the following linear system: Find $u_{hb}^{k+1}\in V_{hb}^{D}$ such that
\begin{equation*}
B(u_{hb}^{k+1},v_{hb})+\sum\limits_{T\in\mathcal{T}_{h}}\int_T\xi_{T}^{k+1}(u_{h}^{k+1})\nabla u_{hb}^{k+1}\cdot \nabla v_{hb}\ d{\bf x}=(f,v_{hb})+(g,v_{hb})_{\Gamma_{N}},\ \forall v_{hb}\in V_{hb}^{D};
\end{equation*}
\item[Step 4] While $k< M$ and
\begin{equation*}
\frac{\|u_{h}^{k+1}-u_{h}^{k}\|_{0,\infty}}{\|u_{h}^{k}\|_{0,\infty}}\geq \delta,
\end{equation*}
set $k=k+1$ and go to Step 2.
\end{itemize}

\subsection{Numerical examples}
In this subsection, we apply the fixed point iteration algorithm described in subsection 6.1 to compute  two numerical  examples to test the convergence rate of the modified DD method  (\ref{formulatwo}). The first example is taken from \cite{Santos2021} and has no   interior or
boundary layer,   and the second one is   from \cite{Du2015}  and has boundary layers when $\varepsilon$ is small.

\begin{figure}[htbp!]
\centering
 {\includegraphics[height=3cm,width=3.5cm]{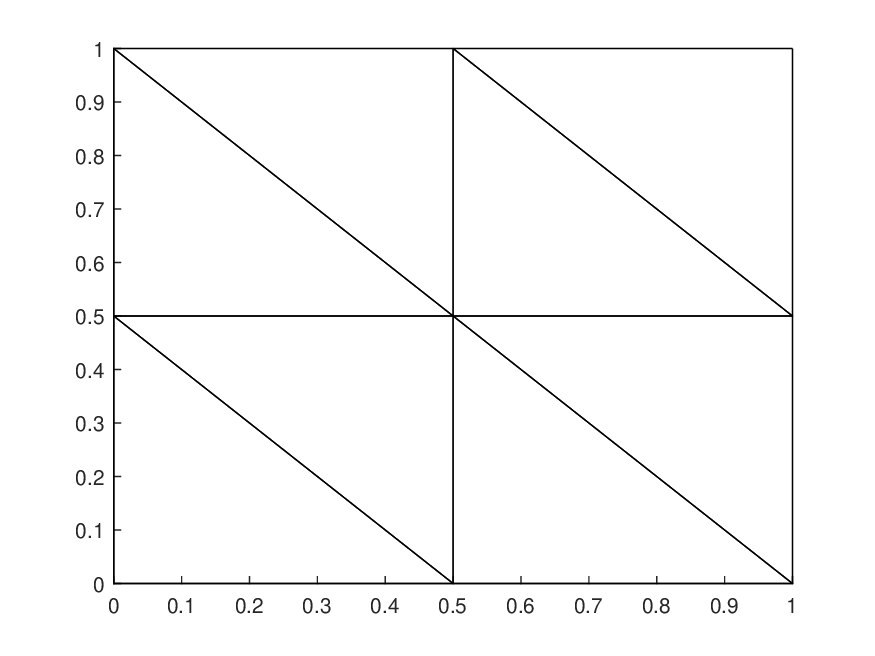}}
 {\includegraphics[height=3cm,width=3.5cm]{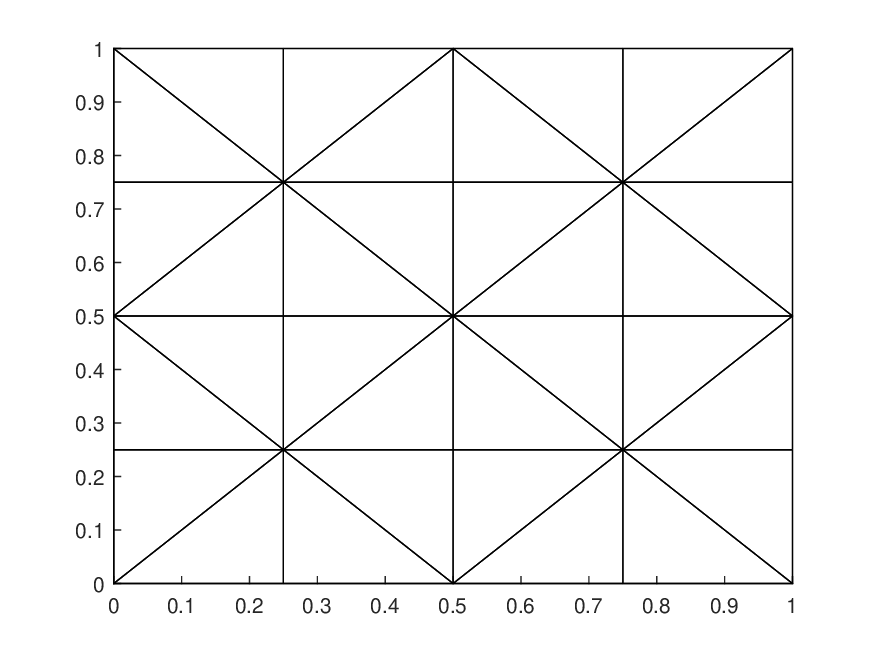}}
{\includegraphics[height=3cm,width=3.5cm]{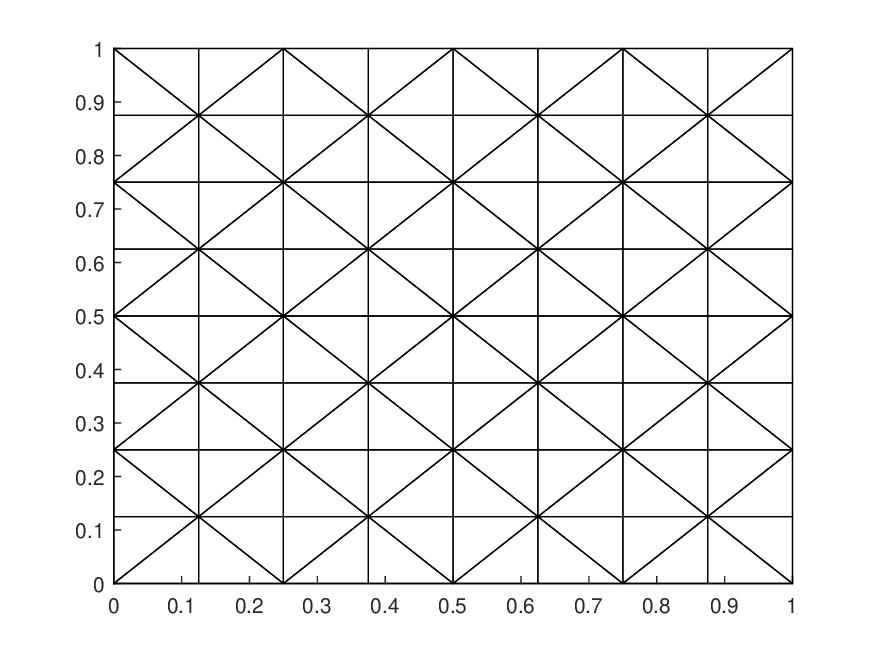}}
\caption{ The domain $\Omega$ and $N\times N$ meshes with   $N=2,\ 4,\ 8$.}
\label{fig61}
\end{figure}

In both of the two examples,   we take $\Omega=(0,1)\times(0,1)$ in the model (\ref{formula1}) with $N\times N$ uniform triangular meshes (cf. Figure \ref{fig61}) for $N=2,4, 8, 16,32,64$.

\begin{example}[A  problem without  interior or boundary layers]\label{exam1}
In this example, we set  $\beta=(3,2)^{T}$ and the exact solution  of (\ref{formula1})  is given by
\begin{equation*}
u(x,y)=100x^{2}(1-x)^{2}y(1-y)(1-2y).
\end{equation*}
For   the setting of parameters    $\varepsilon$ and $\sigma$, we consider four cases: (i) $\varepsilon=10, \sigma=1 $; (ii) $ \varepsilon=10,\sigma=0,$; (iii) $ \varepsilon=10^{-6}, \sigma=1$;  (iv) $\varepsilon=10^{-6},\sigma=0$.

Tables \ref{du1} to \ref{du4}  list some numerical  results of  the errors $\|u-u_{hb}\|_{0}$, $|u-u_{hb}|_{1}$ and
\begin{equation*}
|||u-u_{hb}|||_{*}:=\Big\{|||u-u_{hb}|||^{2}+\sum\limits_{T\in\mathcal{T}_{h}}D_{T}(u_{hb};u_{hb},u_{hb})\Big\}^{1/2}.
\end{equation*}

%
\end{example}

\begin{figure}[htbp]
    \begin{minipage}[t]{0.5\linewidth}
        \centering
        \includegraphics[width=2.25in]{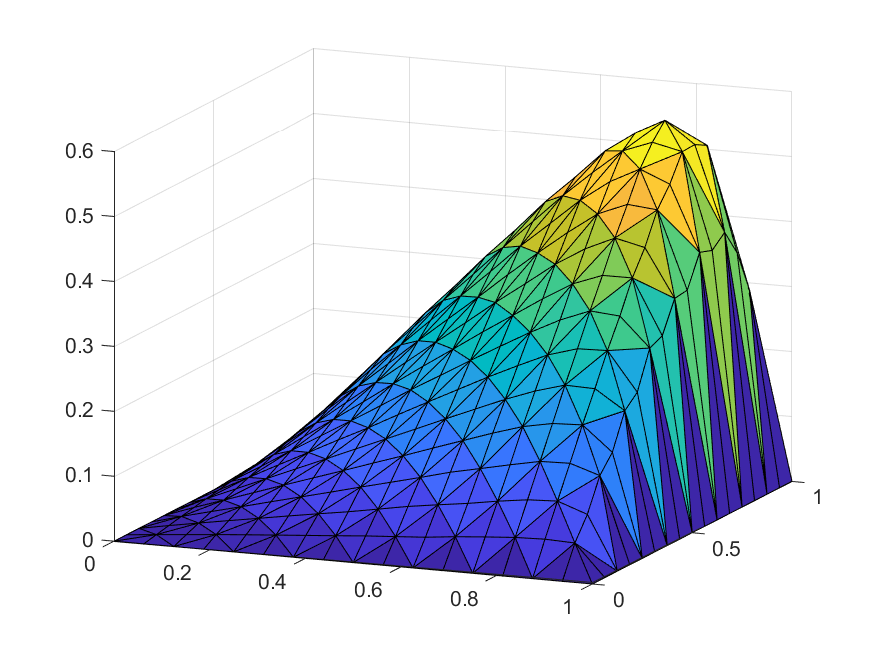}\\
      \end{minipage}
    \begin{minipage}[t]{0.5\linewidth}
        \centering
        \includegraphics[width=2.25in]{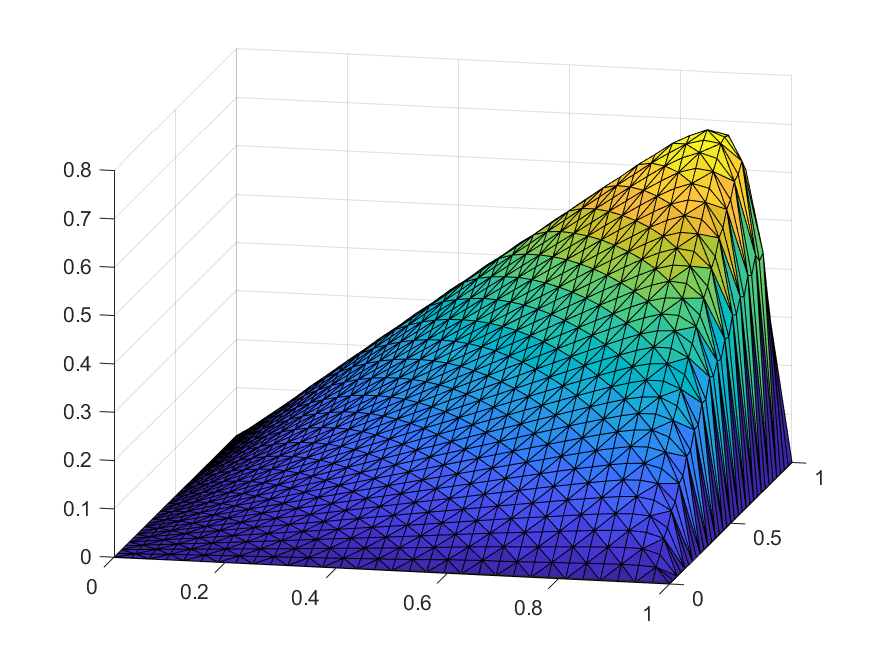}\\
    \end{minipage}
    \caption{Approximation solution   $u_{hb}$ over   meshes $8\times8$ (left) and $16\times16$ (right)
   with $\varepsilon=10^{-6}, \sigma=1$.}
    \label{dutu1}
\end{figure}

\begin{example}[A  problem with    boundary layers]\label{exam2}
We set $\beta=(1,1)^{T}$ and the exact solution is given by
\begin{equation*}
\displaystyle u(x,y)=\left(\frac{\exp\left({\frac{x-1}{\varepsilon}}\right)-1}{\exp\left({-\frac{1}{\varepsilon}}\right)-1}+x-1\right)
\left(\frac{\exp\left({\frac{y-1}{\varepsilon}}\right)-1}{\exp\left({-\frac{1}{\varepsilon}}\right)-1}+y-1\right).
\end{equation*}
we consider three cases of the parameters     $\varepsilon$ and $\sigma$: (i) $\varepsilon=10, \sigma=1 $; (ii) $ \varepsilon=0.1, \sigma=1$; (iii) $ \varepsilon=10^{-6},\sigma=1$. Note that the   solution $u$
has boundary layers of width $O(\varepsilon)$ along $x=1$ and $y=1$.

  Tables \ref{du5} to \ref{du7}   show some numerical  results of  the errors $\|u-u_{hb}\|_{0}$, $|u-u_{hb}|_{1}$ and $|||u-u_{hb}|||_{*}$, and Figure \ref{dutu1} reports the approximation solutions   over the $16\times 16$  and
 $32\times 32$ meshes.

\end{example}

For Example \ref{exam1}, from Tables \ref{du1} to \ref{du4}   we have the following observations:
\begin{itemize}
  \item
Tables \ref{du1} and \ref{du3} show numerical results    in the  diffusion-dominated cases with $\varepsilon=10 $ and $\sigma=1,\ 0$. We can see that  the convergence rates of $\|u-u_{hb}\|_{0}$, $|u-u_{hb}|_{1}$ and $|||u-u_{hb}|||_{*}$ are second, first and first orders, respectively.
%

\item Tables  \ref{du2} and \ref{du4} show numerical results in the convection-dominated
cases with  $\varepsilon=10^{-6}$ and $\sigma=1,\ 0$. We can still see that the convergence rates of $\|u-u_{hb}\|_{0}$, $|u-u_{hb}|_{1}$ and $|||u-u_{hb}|||_{*}$ are   second, first and first   orders, respectively.

\item The convergence rates of the errors $|u-u_{hb}|_{1}$ and $|||u-u_{hb}|||_{*}$ are conformable to the theoretical results in Theorem \ref{formula45}. 

%
\end{itemize}

For Example \ref{exam2}, from Tables \ref{du5} to \ref{du7}    we have the following observations: 
\begin{itemize}
\item Tables \ref{du5} and  \ref{du6} gives numerical results in a diffusion-dominated case with $\sigma=1$ and $\varepsilon=10$ and in a slightly convection$/$reaction dominated case with $\sigma=1$ and $\varepsilon=0.1$.  We can see that the convergence rates of $\|u-u_{hb}\|_{0}$, $|u-u_{hb}|_{1}$ and $|||u-u_{hb}|||_{*}$ are   second, first and first   orders, respectively, as is conformable to the theoretical results in Theorem \ref{formula45}.

\item Table \ref{du7} provides numerical results in a strong convection-dominated case with $\sigma=1$ and $\varepsilon=10^{-6}$. We can see that the three types of errors converge to zero very slowly, due to the boundary layers of  the exact solution. As shown in Remark \ref{pinglun3}, to resolve the layers one may refer to a certain adaptive  finite element method.

\item Figure \ref{dutu1} depicts the approximation solutions over   meshes $8\times8$   (left) and
 $16\times16$(right) in the strong convection-dominated case. We can see that the numerical solutions have almost no oscillations around the boundary
layer regions.

\end{itemize}
\begin{table}[t]\small
 \begin{center}
        \caption{Numerical results of Example \ref{exam1}:  $\varepsilon=10,\sigma=1$}
\label{du1}
        \small 
        \begin{tabular}{|c|c|c|c|c|c|c|} \hline
            mesh& $||u-u_{hb}||$& rate& $|u-u_{hb}|_{1}$& rate& $|||u-u_{hb}|||_{*}$& rate\\ \hline
            $2\times2$&0.3338&--&2.8105&--&8.7438&-- \\ \hline
            $4\times4$&0.0936&1.8344&1.4096&0.9955&4.4488&0.9748\\ \hline
            $8\times8$&0.0210&2.1561&0.6562&1.1031&2.0741&1.1009\\ \hline
            $16\times16$&0.0053&1.9863&0.3292&0.9952&1.0408&0.9948\\ \hline
            $32\times32$&0.0013&2.0275&0.1644&1.0018&0.5198&1.2224\\ \hline
            $64\times64$&0.0003&2.1155&0.0821&1.0018&0.2595&1.2224\\ \hline
        \end{tabular}
    \end{center}
\end{table}

\begin{table}[t]\small
 \begin{center}
        \caption{Numerical results of Example \ref{exam1}: $\varepsilon=10,\sigma=0$} \label{du3}
        \small 
        \begin{tabular}{|c|c|c|c|c|c|c|} \hline
            mesh & $||u-u_{hb}||$& rate& $|u-u_{hb}|_{1}$& rate& $|||u-u_{hb}|||_{*}$& rate\\ \hline
            $2\times2$&0.3617&--&2.8103&--&8.8076&-- \\ \hline
            $4\times4$&0.0937&1.9487&1.4093&0.9911&4.4476&0.9857\\ \hline
            $8\times8$&0.0210&2.1577&0.6561&1.1030&2.0737&1.1008\\ \hline
            $16\times16$&0.0053&1.9863&0.3291&0.9954&1.0406&0.9948\\ \hline
            $32\times32$&0.0013&2.0275&0.1644&1.0013&0.5197&1.0017\\ \hline
            $64\times64$&0.0003&2.1155&0.0820&1.0035&0.2594&1.0025\\ \hline
        \end{tabular}
    \end{center}
\end{table}

\begin{table}[t]\small
 \begin{center}
        \caption{Numerical results of Example \ref{exam1}:  $\varepsilon=10^{-6}, \sigma=1$} \label{du2}
        \small 
        \begin{tabular}{|c|c|c|c|c|c|c|} \hline
            mesh & $||u-u_{hb}||$& rate& $|u-u_{hb}|_{1}$&rate& $|||u-u_{hb}|||_{*}$& rate\\ \hline
            $2\times2$&0.3211&--&3.8912&--&1.6094&-- \\ \hline
            $4\times4$&0.0856&1.9073&1.9915&0.9664&0.8184&0.9756\\ \hline
            $8\times8$&0.0195&2.1341&0.9815&1.0208&0.4006&1.0306\\ \hline
            $16\times16$&0.0045&2.1155&0.4912&0.9987&0.2059&0.9602\\ \hline
            $32\times32$&0.0011&2.0324&0.2452&1.0024&0.1039&0.9867\\ \hline
            $64\times64$&0.0003&2.0265&0.1225&1.0012&0.0522&0.9931\\ \hline
        \end{tabular}
    \end{center}
\end{table}

\begin{table}[t]\small
 \begin{center}
        \caption{Numerical results of Example \ref{exam1}:  $\varepsilon=10^{-6},\sigma=0$} \label{du4}
        \small 
        \begin{tabular}{|c|c|c|c|c|c|c|} \hline
            mesh& $||u-u_{hb}||$& rate& $|u-u_{hb}|_{1}$& rate& $|||u-u_{hb}|||_{*}$& rate\\ \hline
            $2\times2$&0.3210&--&3.9012&--&1.6080&-- \\ \hline
            $4\times4$&0.0839&1.9358&1.9693&0.9862&0.8114&0.9868\\ \hline
            $8\times8$&0.0199&2.0759&0.9833&1.0020&0.4024&1.0118\\ \hline
            $16\times16$&0.0046&2.1131&0.4862&1.0161&0.2013&0.9993\\ \hline
            $32\times32$&0.0011&2.0641&0.2427&1.0024&0.1006&1.0007\\ \hline
            $64\times64$&0.0003&2.0265&0.1212&1.0018&0.0501&0.9971\\ \hline
        \end{tabular}
    \end{center}
\end{table}

\begin{table}[t]\small
 \begin{center}
        \caption{Numerical results of Example \ref{exam2}: $\varepsilon=10,\sigma=1$ } \label{du5}
        \small 
        \begin{tabular}{|c|c|c|c|c|c|c|} \hline
            mesh& $||u-u_{hb}||$& rate& $|u-u_{hb}|_{1}$& rate& $|||u-u_{hb}|||_{*}$& rate\\ \hline
            $2\times2$&3.2e-05&--&2.7e-04&--&8.5e-04&-- \\ \hline
            $4\times4$&9.9e-06&1.6926&1.3e-04&0.9952&4.3e-04&0.9956\\ \hline
           $8\times8$&2.6e-06&1.9289&7.0e-05&0.9494&2.2e-04&0.9494\\ \hline
           $16\times16$&6.6e-07&1.9780&3.5e-05&0.9938&1.1e-04&0.9928\\ \hline
            $32\times32$&1.6e-07&2.0444&1.8e-05&1.0000&5.5e-05&1.0000\\ \hline
            $64\times64$&4.0e-08&2.0000&8.8e-06&0.9918&2.8e-05&1.0000\\ \hline
        \end{tabular}
    \end{center}
\end{table}

\begin{table}[t]\small
 \begin{center}
        \caption{Numerical results of Example \ref{exam2}: $\varepsilon=0.1,\sigma=1$ } \label{du6}
        \small 
        \begin{tabular}{|c|c|c|c|c|c|c|} \hline
            mesh& $||u-u_{hb}||$& rate& $|u-u_{hb}|_{1}$&rate& $|||u-u_{hb}|||_{*}$& rate \\ \hline
            $2\times2$&0.1644&--&1.9810&--&0.6049&-- \\ \hline
           $4\times4$&0.0417&1.9791&1.0005&0.9855&0.3201&0.9182\\ \hline
            $8\times8$&0.0120&1.7970&0.5005&0.9993&0.1612&0.9897\\ \hline
            $16\times16$&0.0034&1.8194&0.2676&0.9033&0.0847&0.9284\\ \hline
            $32\times32$&0.0009&1.9175&0.1343&0.9946&0.0425&0.9949\\ \hline
            $64\times64$&0.0002&2.1699&0.0668&1.0075&0.0211&1.0102\\ \hline
        \end{tabular}
    \end{center}
\end{table}

\begin{table}[t]\small
 \begin{center}
        \caption{Numerical results of Example \ref{exam2}: $\varepsilon=10^{-6}, \sigma=1 $ } \label{du7}
        \small 
        \begin{tabular}{|c|c|c|c|c|} \hline
            mesh &$||u-u_{hb}||_{0}$&rate &$|u-u_{hb}|_{1}$&$|||u-u_{hb}|||_{*}$\\ \hline
           $2\times2$&0.3049&--&3.2874&0.5847\\ \hline
            $4\times4$&0.2593&0.2331&2.8286&0.5550\\ \hline
            $8\times8$&0.2049&0.3397&1.9776&0.5203\\ \hline
            $16\times16$&0.1548&0.4045&1.3927&0.4694\\ \hline
           $32\times32$&0.1136&0.4464&1.0316&0.4646\\ \hline
            $64\times64$&0.0809&0.4898&0.9336&0.4119\\ \hline
        \end{tabular}
    \end{center}
\end{table}

\end{document}